\documentclass[12pt]{article}
\usepackage{amsmath}
\usepackage{amssymb}
\usepackage{amscd}
\usepackage{latexsym}
\usepackage{theorem}
\usepackage{doublespace}
\usepackage{epsfig}
\usepackage{amsfonts}
\usepackage{enumerate}

\newtheorem{theorem}{Theorem}[subsection]

\newtheorem{proposition}[theorem]{Proposition}
\newtheorem{corollary}[theorem]{Corollary}

\newtheorem{conjecture}[theorem]{Conjecture}

{
\theorembodyfont{\rmfamily}
\theoremstyle{plain}

\newtheorem{example}[theorem]{Example}

\newtheorem{remark}[theorem]{Remark}
}

\renewcommand{\Re}{\operatorname{Re}}
\renewcommand{\Im}{\operatorname{Im}}

\renewcommand{\dim}{\operatorname{dim}}

\newcommand{\Sym}{\operatorname{Sym}}

\newcommand{\C}{{\mathbb{C}}}
\newcommand{\Z}{{\mathbb{Z}}}

\newcommand{\R}{{\mathbb{R}}}

\newcommand{\otn}{\{1,\ldots,n\}}

\newcommand{\bigmid}{\hs\Big{|}\hs}
\newcommand{\subs}{\subseteq}

\newcommand{\hookto}{{\hookrightarrow}}

\newcommand{\hs}{\hspace{3pt}}

\newcommand{\D}{\Delta}
\renewcommand{\i}{\iota}
\renewcommand{\a}{\alpha}

\renewcommand{\cot}{T^*\Cn}

\newcommand{\Cn}{\C^n}

\newcommand{\M}{\mathfrak{M}}

\newcommand{\X}{\mathfrak{X}}

\newcommand{\minv}{\mu^{\! -1}(0)}

\renewcommand{\mod}{{/\!\!/}}
\newcommand{\mmod}{{\!/\!\!/\!\!/\!\!/}}

\newcommand{\gd}{\g^*}
\newcommand{\g}{\mathfrak{g}}

\newcommand{\Tk}{T^k}
\newcommand{\Tn}{T^n}
\newcommand{\Td}{T^d}
\newcommand{\tk}{\mathfrak{t}^k}
\newcommand{\tn}{\mathfrak{t}^n}
\newcommand{\td}{\mathfrak{t}^d}
\newcommand{\tkd}{(\tk)^*}

\newcommand{\tdz}{\td_{\Z}}
\newcommand{\tnd}{(\tn)^*}

\newcommand{\tddr}{(\td)_{\R}^*}
\newcommand{\tddz}{(\td)_{\Z}^*}

\newcommand{\A}{\mathcal{A}}

\newcommand{\At}{\tilde\A}

\newcommand{\MAt}{\M(\At)}

\newcommand{\MA}{\M(\A)}

\newcommand{\da}{\Delta_\A}
\newcommand{\bcsda}{\operatorname{bc}_\sigma\!\da}

\newcommand{\Proj}{\mathrm{Proj}\,}

\newcommand{\IH}{I\! H}

\newcommand{\IC}{\mathbf{IC}}

\newcommand{\DA}{\Delta_{\A}}

\newcommand{\SR}{\mathcal{SR}}

\newcommand{\Mal}{\M_{\a,\la}}
\newcommand{\la}{\lambda}
\newcommand{\Ml}{\M_{\la}}

\newcommand{\<}{\left<}
\renewcommand{\>}{\right>}

\newcommand{\Poin}{\operatorname{Poin}}

\newcommand{\vlf}{(\Cn)^{\ell\! f}}
\newcommand{\xlf}{\X^{\ell\! f}}
\newcommand{\vast}{(\Cn)^{\a-st}}
\newcommand{\vass}{(\Cn)^{\a-ss}}

\newcommand{\tkdz}{\tkd_\Z}
\newcommand{\cs}{\C^\times}
\newcommand{\surj}{\twoheadrightarrow}
\newcommand{\Maz}{\M_{\a,0}}
\newcommand{\Mzz}{\M_{0,0}}
\renewcommand{\t}{\mathfrak t}

\newcommand{\qed}{\hfill \mbox{$\Box$}\medskip\newline}

\newenvironment{proofsmoothring}{\noindent {\bf Proof of \ref{smoothring}:}}{\qed \par}

\begin{document}
\spacing{1.1}
\noindent
{\Large \bf A survey of hypertoric geometry and topology}
\bigskip\\
{\bf Nicholas Proudfoot}\footnote{Supported
by a National Science Foundation Mathematical Sciences Postdoctoral Research Fellowship.}\\
Department of Mathematics, Columbia University, 10027
\bigskip
{\small
\begin{quote}
\noindent {\em Abstract.}
Hypertoric varieties are quaternionic analogues of toric varieties, important
for their interaction with the combinatorics of matroids as well as for their
prominent place in the rapidly expanding field of algebraic symplectic and
hyperk\"ahler geometry.  The aim of this survey is to give clear definitions and statements
of known results, serving both as a reference and as a point of entry to this beautiful subject.
\end{quote}
}
\bigskip

\noindent
Given a linear representation of a reductive complex algebraic group $G$,
there are two natural quotient constructions.  First, one can take a geometric
invariant theory (GIT) quotient, which may also be interpreted as a K\"ahler
quotient by a maximal compact subgroup of $G$.  Examples of this sort
include toric varieties (when $G$ is abelian), moduli spaces of spacial polygons,
and, more generally, moduli spaces of semistable representations of quivers.
A second construction involves taking an algebraic symplectic quotient of
the cotangent bundle of $V$, which may also be interpreted as a hyperk\"ahler
quotient.  The analogous examples of the second type are hypertoric varieties,
hyperpolygon spaces, and Nakajima quiver varieties.

The subject of this survey will be hypertoric varieties, which are by definition
the varieties obtained from the second construction when $G$ is abelian.
Just as the geometry and topology of toric varieties is deeply connected
to the combinatorics of polytopes, hypertoric varieties interact richly with
the combinatorics of hyperplane arrangements and matroids.  Furthermore,
just as in the toric case, the flow of information goes in both directions.

On one hand, Betti numbers of hypertoric varieties have a combinatorial 
interpretation, and the geometry of the varieties can be used to prove
combinatorial results.  Many purely algebraic constructions involving matroids
acquire geometric meaning via hypertoric varieties, and this has led to geometric proofs of
special cases of the g-theorem for matroids \cite[7.4]{HSt} and the Kook-Reiner-Stanton
convolution formula \cite[5.4]{PW}.  Future plans include a geometric interpretation
of the Tutte polynomial and of the phenomenon of Gale duality of matroids \cite{BLP}.

On the other hand, hypertoric varieties are important to geometers with no interest
in combinatorics simply because they are among the most explicitly understood
examples of algebraic symplectic or hyperk\"ahler varieties, which are becoming
increasingly prevalent in many areas of mathematics.  For example,
Nakajima's quiver varieties include resolutions of Slodowy slices and 
Hilbert schemes of points on ALE spaces, both of which play major roles
in modern representation theory.  Moduli spaces of Higgs bundles are currently
receiving a lot of attention in string theory, and character varieties of fundamental
groups of surfaces and 3-manifolds have become an important tool in low-dimensional
topology.  Hypertoric varieties are useful for understanding such spaces
partly because their geometries share various characteristics, and partly
through explicit abelianization theorems, such as those stated and conjectured
in Section \ref{abel}.

Our main goal is to give clear statements of the definitions and selected theorems that 
already appear in the literature, along with explicit references.
With the exception of Theorem \ref{smoothring}, for which we give a new proof,
this article does not contain any proofs at all.
Section \ref{ht} covers the definition of hypertoric varieties, and explains
their relationship to hyperplane arrangements.  Section \ref{homotopy}
gives three different constructions of unions of toric varieties that are equivariantly
homotopy equivalent to a given hypertoric variety.  These constructions have
been extremely useful from the standpoint of computing algebraic invariants,
and can also make hypertoric varieties more accessible to someone with a
background in toric geometry but less experience with algebraic symplectic
or hyperk\"ahler quotients.  Finally, Section \ref{cohomology} in concerned
with the cohomology of hypertoric varieties, giving concrete form to the general principle
that hypertoric geometry is intricately related to the combinatorics of matroids.

Section \ref{homotopy} assumes a familiarity with toric varieties, but Sections \ref{ht}
and \ref{cohomology} can both be read independently of Section \ref{homotopy}.
The main quotient construction of Section \ref{hkas}
is logically self-contained, but may be fairly opaque
to a reader who is not familiar with geometric invariant theory.  
Two alternative interpretations of this construction are given in
Remarks \ref{hk} and \ref{cover}, or one can take it as a black box
and still get a sense of the combinatorial flavor of the subject.

\paragraph{\bf Acknowledgments.}  The author is grateful to the organizers of the 2006
International Conference on Toric Topology at Osaka City University, out of which this survey grew.

\begin{section}{Definitions and basic properties}\label{ht}
Hypertoric varieties can be considered either as algebraic varieties or, in the smooth case, 
as hyperk\"ahler manifolds.  In this section we give a constructive definition, with a strong
bias toward the algebraic interpretation.  Section \ref{hkas} proceeds in greater generality
than is necessary for hypertoric varieties so as to unify the theory with that of other algebraic
symplectic quotients, most notably Nakajima quiver varieties.

\begin{subsection}{Algebraic symplectic quotients}\label{hkas}
Let $G$ be a reductive algebraic group over the complex numbers 
acting linearly and effectively
on a finite-dimensional complex vector space $V$.
The cotangent bundle $$T^*V \cong V \times V^*$$ carries a natural 
algebraic symplectic form $\Omega$.
The induced action of $G$ on $T^*V$ is hamiltonian, with moment map
$$\mu:T^*V\to\gd$$ given by the equation
$$\mu(z,w)(x) =  \Omega(x\cdot z, w)\,\,
\text{ for all $z\in V$, $w\in V^*$, $x\in\g$.}$$
Suppose given an element $\la\in Z(\gd)$ (the part of $\gd$ fixed by
the coadjoint action of $G$), and a multiplicative character $\a:G\to\cs$,
which may be identified with an element of $Z(\gd_\Z)$ by taking its derivative
at the identity element of $G$.\footnote{Strictly speaking, an element
of $Z(\gd_\Z)$ only determines a character of the connected component
of the identity of $G$.  It can be checked, however, that the notion of $\a$-stability
defined below depends only on the restriction of $\a$ to the identity component,
therefore we will abusively think of $\a$ as sitting inside of $Z(\gd_\Z)$.}
The fact that $\la$ lies in $Z(\gd)$ implies
that $G$ acts on $\mu^{-1}(\la)$.
Our main object of study in this survey will be
the algebraic symplectic quotient 
\begin{equation*}
\Mal=T^*V\mmod_{\!\!\a}G = \mu^{-1}(\la)\mod_{\!\!\a}G.
\end{equation*}
Here the second quotient is a projective GIT quotient
\begin{equation}\label{gitdef}\mu^{-1}(\la)\mod_{\!\!\a}G := \Proj\bigoplus_{m=0}^{\infty}
\Big\{f\in\operatorname{Fun}\big(\mu^{-1}(\la)\big)
\bigmid \nu(g)^*f = \a(g)^m f\hs\text{ for all }g\in G\hs\Big\},\end{equation}
where $\nu(g)$ is the automorphism of $\mu^{-1}(\la)$ defined by $g$.

This quotient may be defined in a more geometric way as follows.  A point
$(z,w)\in\mu^{-1}(\la)$ is called {\bf $\mathbf\a$-semistable} if there exists a function
$f$ on $\mu^{-1}(\la)$ and a positive integer $m$ such that $\nu(g)^*f = \a(g)^m f$
for all $g\in G$ and $f(z,w)\neq 0$.  
It is called {\bf $\mathbf\a$-stable} if it is $\a$-semistable and its
$G$-orbit in the $\a$-semistable set is closed with finite stabilizers.  Then
the stable and semistable sets
$$\mu^{-1}(\la)^{\a-st}\subs\mu^{-1}(\la)^{\a-ss}\subs\mu^{-1}(\la)$$
are nonempty and Zariski open, and there is a surjection
$$\mu^{-1}(\la)^{\a-ss}\surj\Mal$$ with $(z,w)$ and $(z',w')$ mapping to the
same point if and only if the closures of their $G$-orbits intersect in $\mu^{-1}(\la)^{\a-ss}$.
In particular, the restriction of this map to the stable locus is nothing but the geometric
quotient by $G$.  For an introduction to geometric invariant theory that explains
the equivalence of these two perspectives, see \cite{P2}.

\begin{remark}\label{hk}
The algebraic symplectic quotient defined above may also be interpreted
as a hyperk\"ahler quotient.  The even dimensional complex vector space
$T^*V$ admits a complete hyperk\"ahler metric, and the action of the maximal
compact subgroup $G_\R\subs G$ is {\bf hyperhamiltonian}, meaning that it is
hamiltonian with respect to all three of the real symplectic forms on $T^*V$.
Then $\Mal$ is naturally diffeomorphic to the hyperk\"ahler quotient of $T^*V$ by $G_\R$,
in the sense of \cite{HKLR}, at the value $(\a,\Re\la,\Im\la)\in\gd_\R\otimes\R^3$.
This was the original perspective on both hypertoric varieties \cite{BD} and
Nakajima quiver varieties \cite{N1}.  For more on this perspective in the hypertoric case,
see Konno's survey in this volume \cite[\S 3]{K4}.
\end{remark}

We note that if $\a = 0$ is the trivial character of $G$, then Equation \eqref{gitdef}
simplifies to $$\M_{0,\la} = \operatorname{Spec}\operatorname{Fun}
\big(\mu^{-1}(\la)\big)^{G}.$$
Furthermore, since $\Mal$ is defined as the projective spectrum of a graded ring whose
degree zero part is the ring of invariant functions on $\mu^{-1}(\la)$, we always have a projective morphism
\begin{equation}\label{projmorph}\Mal\surj\M_{0,\la}.\end{equation}
This morphism may also be induced from the inclusion of the inclusion
$$\mu^{-1}(\la)^{\a-ss}\subs\mu^{-1}(\la) = \mu^{-1}(\la)^{0-ss}.$$
From this we may conclude that it is generically one-to-one, and therefore a partial resolution.
When $\la=0$, we have a distinguished point in $\M_{0,0}$, namely the image of $0\in\minv$
under the map induced by the inclusion of the invariant functions into the
coordinate ring of $\minv$.  The preimage of this point under
the morphism \eqref{projmorph} is called the {\bf core} of $\M_{\a,0}$, and will be further
studied (in the case where $G$ is abelian) in Section \ref{core}.

On the other extreme, if $\la$ is a regular value of $\mu$, then $G$ will act
locally freely on $\mu^{-1}(\la)$.  In this case {\em all} points will be $\a$-stable
for any choice of $\a$, and the GIT quotient
\begin{equation*}\label{ml}\Ml = \mu^{-1}(\la)\mod G\end{equation*}
will simply be a geometric quotient.  In particular, the morphism
\eqref{projmorph} becomes an isomorphism.
Both the case of regular $\la$ and the case $\la=0$ will be of interest to us.

We call a pair $(\a,\la)$ {\bf generic} if $\mu^{-1}(\la)^{\a-st} = \mu^{-1}(\la)^{\a-ss}$.
In this case the moment map condition tells us that the stable set is smooth, and therefore
that the quotient $\Mal$ by the locally free $G$-action has at worst orbifold singularities.
Using the hyperk\"ahler quotient perspective of Remark \ref{hk}, one can prove the
following Proposition.  (See \cite[2.6]{K3} or \cite[2.1]{HP1} in the hypertoric case, and
\cite[4.2]{N1} in the case of quiver varieties; the general case is no harder than these.)

\begin{proposition}\label{diffeo}
If $(\a,\la)$ and $(\a',\la')$ are both generic, then 
the two symplectic quotients $\Mal$
and $\M_{\a',\la'}$ are diffeomorphic.
\end{proposition}

\begin{remark}  If $G$ is semisimple, 
then $Z(\gd) = \{0\}$, and (unless $G$ is finite) it will not be possible
to choose a regular value $\la\in Z(\gd)$, nor a nontrivial character $\a$.  
We will very soon specialize, however, 
to the case where $G$ is abelian.  In this case $Z(\gd) = \gd$, the regular
values form a dense open set, and the characters of $G$ form a full integral lattice 
$\gd_\Z\subs\gd$.
\end{remark}
\end{subsection}

\begin{subsection}{Hypertoric varieties defined}\label{defs}
Let $\tn$ be the coordinate complex vector space
of dimension $n$ with basis $\{\varepsilon_1,\ldots,\varepsilon_n\}$, and let $\td$
be a complex vector space of dimension $d$ with a full lattice $\tdz$.
Though $\tdz$ is isomorphic to the standard integer lattice $\Z^d$,
we will not choose such an isomorphism.
Let $\{a_1,\ldots,a_n\}\subset\tdz$ be a collection of nonzero vectors
such that the map $\tn\to\td$ taking $\varepsilon_i$ to $a_i$ is surjective.  
Let $k=n-d$, and let $\tk$ be the kernel of this map.
Then we have an exact sequence
\begin{eqnarray}\label{vs}
0 \longrightarrow \tk \stackrel{\i}{\longrightarrow} \tn
\longrightarrow \td\longrightarrow 0,
\end{eqnarray}
which exponentiates to an exact sequence of tori
\begin{eqnarray}\label{tori}
0 \longrightarrow \Tk \longrightarrow \Tn
\longrightarrow \Td\longrightarrow 0.
\end{eqnarray}
Here $\Tn=(\cs)^n$, $\Td$ is a quotient of $\Tn$, and
$\Tk=\ker\!\big(\Tn\to\Td\big)$ is a subgroup with
Lie algebra $\tk$, which is connected if and only if the vectors
$\{a_i\}$ span the lattice $\tdz$ over the integers.  Note that every algebraic
subgroup of $\Tn$ arises in this way.

The torus $\Tn$ acts naturally via coordinatewise multiplication on the vector space $\Cn$,
thus so does the subtorus $\Tk$.
For $\a\in\tkdz$ a multiplicative character of $\Tk$ and $\la\in\tkd$ arbitrary, the  
algebraic symplectic quotient
\begin{equation*}
\Mal = \cot\mmod_{\!\!\a}\Tk
\end{equation*}
is called a {\bf hypertoric variety}.

The hypertoric variety $\Mal$ is a symplectic variety of dimension $2d$ which admits
a complete hyperk\"ahler metric.  The action of the quotient torus $\Td=\Tn/\Tk$ on $\Mal$
is hamiltonian with respect to the algebraic symplectic form, and the action of the maximal
compact subtorus $\Td_\R$ is hyperhamiltonian.
In the original paper of Bielawski and Dancer \cite{BD} the hyperk\"ahler perspective
was stressed, and the spaces were referred to as ``toric hyperk\"ahler manifolds".
However, since we have worked frequently with singular reductions as well as with fields of definition
other than the complex numbers (see for example \cite{HP1,P3,PW}), we prefer the term
hypertoric varieties.

\begin{remark}\label{eq}
In the hypertoric case, the diffeomorphism of Proposition \ref{diffeo} 
can be made $\Td_\R$-equivariant \cite[2.1]{HP1}.
\end{remark}
\end{subsection}

\begin{subsection}{Hyperplane arrangements}\label{arrangements}
The case in which $\la=0$ will be of particular importance, and 
it is convenient to encode the data that were used to construct the hypertoric variety
$\M_{\a,0}$ in terms of an arrangement of affine hyperplanes 
with some additional structure in the real vector space
$\tddr = \tddz\otimes_\Z\R$.  A {\bf weighted,
cooriented, affine hyperplane} $H\subs\tddr$ is an affine
hyperplane along with a choice of nonzero integer normal vector
$a\in\tdz$.  Here ``affine'' means that $H$ need not pass through the origin, and
``weighted'' means that $a$ is not required to be primitive. Let
$r = (r_1,\ldots,r_n)\in\tnd$ be a lift of $\a$ along $\i^*$, and
let $$H_i =\{x\in\tddr \mid x\cdot a_i + r_i = 0\}$$ be the
weighted, cooriented, affine hyperplane with normal vector
$a_i\in\tddz$.  (Choosing a different $r$ corresponds to simultaneously
translating all of the hyperplanes by a vector in $\tddz$.)
We will denote the collection $\{H_1,\ldots,H_n\}$
by $\A$, and write
$$\MA = \M_{\a,0}$$ for the corresponding hypertoric variety.
We will refer to $\A$ simply as an {\bf arrangement}, always assuming that the
weighted coorientations are part of the data.

\begin{remark}
We note that we allow repetitions of hyperplanes in our
arrangement ($\A$ may be a multi-set), and that a repeated
occurrence of a particular hyperplane is {\em not} the same as a
single occurrence of that hyperplane with weight 2. On the other
hand, little is lost by restricting one's attention to
arrangements of distinct hyperplanes of weight one.
\end{remark}

Since each hyperplane $H_i$ comes with a normal vector, it seems at first
that it would make the most sense to talk about an arrangement of half-spaces,
where the $i^\text{th}$ half-space consists of the set of points that lie on the positive
side of $H_i$ with respect to $a_i$.  The reason that we talk about hyperplanes
rather than half-spaces is the following proposition, proven in \cite[2.2]{HP1}.

\begin{proposition}\label{flip}
The $\Td$-variety $\MA$ does not depend on the signs of the vectors $a_i$.
\end{proposition}

In other words, if we make a new hypertoric variety with the same arrangement of 
weighted hyperplanes but with some of the coorientations flipped, it will be $\Td$-equivariantly
isomorphic to the hypertoric variety with which we started.\footnote{In \cite{HP1} we consider
an extra $\cs$ action on $\MA$ that {\em does} depend on the coorientations.}

We call the arrangement $\A$ {\bf simple} if every subset of $m$
hyperplanes with nonempty intersection intersects in codimension
$m$. We call $\A$ {\bf unimodular} if
every collection of $d$ linearly independent vectors
$\{a_{i_1},\ldots,a_{i_d}\}$ spans $\td$ over the integers.
An arrangement which is both simple and unimodular is called {\bf smooth}.
The following proposition is proven in \cite[3.2 $\&$ 3.3]{BD}.

\begin{proposition}\label{simple}
The hypertoric variety $\MA$ has at worst orbifold (finite quotient)
singularities if and only if $\A$ is simple, and is smooth if and only if $\A$ is smooth.
\end{proposition}

For the remainder of the paper, 
Let $\A=\{H_1,\ldots,H_n\}$ be a {\bf central} arrangement, meaning that $r_i = 0$ for all $i$,
so that all of the hyperplanes pass through the origin.
Then $\MA$ is the singular affine variety $\M_{0,0}$.
Let $\At = \{\tilde H_1,\ldots, \tilde H_n\}$ be a {\bf simplification} of $\A$, by
which we mean an arrangement defined by the same vectors
$\{a_i\}\subset\td$, but with a different choice of $r\in\tnd$, such that
$\At$ is simple. This corresponds to translating each of the
hyperplanes in $\A$ away from the origin by some generic amount.
Then $\MAt$ maps $T$-equivariantly to $\MA$ by Equation \eqref{projmorph}, and
Proposition \ref{simple} tell us that it is in fact an ``orbifold resolution", 
meaning a projective morphism, generically one-to-one, in which the source has at worst
orbifold singularities.  The structure of this map is studied extensively in \cite{PW}.
\end{subsection}

\begin{subsection}{Toward an abstract definition}\label{abstract}
The definition of a hypertoric variety in Section \ref{defs} is constructive, modeled on
the definition of toric varieties as GIT quotients of the form $\Cn\mod_{\!\!\a}\Tk$, or equivalently
as symplectic quotients by compact tori.  In the case of toric varieties, there are also abstract
definitions.  
In the symplectic world, one defines a toric orbifold to be a symplectic orbifold of dimension $2d$ 
along with an effective Hamiltonian action of a compact $d$-torus, and proves that any
connected, compact toric orbifold arises from the symplectic quotient construction \cite{De,LT}.
In the algebraic world, one defines a toric variety to be a normal variety
admitting a torus action with a dense orbit, and then proves that any 
semiprojective\footnote{Hausel and Sturmfels call a toric variety semiprojective
if it is projective over its affinization and has at least one torus fixed point.}
toric variety with at worst orbifold singularities arises from the GIT construction.
This idea goes back to \cite{Co}, and can be found in this language in \cite[2.6]{HSt}.

It is natural to ask for such an abstract definition and classification theorem for hypertoric
varieties, either from the standpoint of symplectic algebraic geometry or that of hyperk\"ahler
geometry.  In the hyperk\"ahler setting, such a theorem was proven in \cite[3,4]{Bi}.

\begin{theorem}\label{bielawski}
Any complete, connected, hyperk\"ahler manifold of real dimension $4d$
which admits an effective, hyperhamiltonian action of the compact torus $\Td_\R$
is $\Td_\R$-equivariantly diffeomorphic, and Taub-NUT deformation equivalent, 
to a hypertoric variety.  Any such manifold with Euclidean volume growth is $\Td_\R$-equivariantly
isometric to a hypertoric variety.
\end{theorem}

An analogous algebraic theorem has not been proven, but it should look something like the
following.

\begin{conjecture}
Any connected, symplectic, algebraic variety which is projective over its affinization
and admits an effective, hamiltonian
action of the algebraic torus $\Td$ is equivariantly isomorphic to a Zariski open subset
of a hypertoric variety.
\end{conjecture}
\end{subsection}
\end{section}

\begin{section}{Homotopy models}\label{homotopy}
In this section we fix the vector configuration $\{a_1,\ldots a_n\}\subs\tdz$,
consider three spaces that are $\Td$-equivariantly homotopy
equivalent to the hypertoric variety $\Mal$ for generic choice of $(\a,\la)$.
Each space is essentially toric rather than hypertoric in nature, and therefore
may provide a way to think about hypertoric varieties in terms of more familiar objects.
Recall that if $\la=0$ then $\Mal = \MAt$ for a simple hyperplane
arrangement $\At$, in which the positions of the hyperplanes (up to simultaneous translation)
are determined by $\a$.  If, on the other hand, $\la$ is a regular value, then
$\Mal = \Ml$ is independent of $\a$.

\begin{subsection}{The core}\label{core}
Recall from Section \ref{arrangements} that we have an equivariant 
orbifold resolution
$$\MAt\to\MA,$$ and from Section \ref{hkas} that the fiber $\mathfrak L(\At)\subs\MAt$
over the most singular point of $\MA$ is called the {\bf core} of $\MAt$.  
The primary interest in the core
comes from the following proposition, originally proven in \cite[6.5]{BD} from the
perspective of Proposition \ref{wholecore}.

\begin{proposition}\label{coreretract}
The core $\mathfrak{L}(\At)$ is a $\Td_\R$-equivariant deformation retract of $\MAt$.
\end{proposition}

\begin{remark}
In fact, Proposition \ref{coreretract} holds in the greater generality of Section \ref{hkas},
for algebraic symplectic quotients $\M_{\a,0}$ by arbitrary reductive groups \cite[2.8]{P1}.
The cores of Nakajima's quiver varieties play an important role in representation theory,
because the fundamental classes of the irreducible components form a natural
basis for the top nonvanishing homology group of $\M_{\a,0}$, which may be interpreted
as a weight space of an irreducible representation of a Kac-Moody algebra \cite[10.2]{N2}.
\end{remark}

We now give a toric interpretation of $\mathfrak L(\At)$.
For any subset $U\subs\otn$,
let 
\begin{equation}\label{pu}
P_U = \{x\in\tddr\mid x\cdot a_i + r_i \geq 0\text{ if $i\in U$ and }
x\cdot a_i + r_i \leq 0\text{ if $i\notin U$}\}.
\end{equation}  
Thus $P_U$ is the polyhedron ``cut out"
by the cooriented hyperplanes of $\At$ after reversing the coorientations of the hyperplanes
with indices in $U$.  Since $\At$ is a weighted arrangement, $P_U$ is a labeled polytope
in the sense of \cite{LT}.  Let 
$$\mathcal E_U = \{(z,w)\in\cot\mid w_i = 0\text{ if $i\in U$ and }z_i=0\text{ if $i\notin U$}\}$$
and $$\mathfrak X_U = \mathcal E_U\mod_{\!\!\a}\Tk.$$
Then $\mathcal E_U\subs\mu^{-1}(0)$,
and therefore $$\mathfrak X_U = \mathcal E_U\mod_{\!\!\a}\Tk\subs\mu^{-1}(0)\mod_{\!\!\a}\Tk=\MAt.$$
The following proposition is proven in \cite[6.5]{BD}, but is stated more explicitly
in this language in \cite[3.8]{P1}.

\begin{proposition}\label{corecomps}
The variety $\mathfrak X_U$ is isomorphic to the toric orbifold classified
by the weighted polytope $P_U$.
\end{proposition}

It is not hard to see that the subvariety $\mathcal E_U\mod_{\!\!\a}\Tk\subs\MAt$ lies inside
the core $\mathfrak L(\At)$ of $\MAt$.  In fact, these subvarieties make up the entire core, as
can be deduced from \cite[\S 6]{BD}.

\begin{proposition}\label{wholecore}
\,\,\,$\mathfrak L(\At) \,\,\,= \displaystyle\bigcup_{P_U\text{ bounded}} \mathcal E_U\mod_{\!\!\a}\Tk
\,\,\subs\,\,\MAt.$
\end{proposition}

Thus $\mathfrak L(\At)$ is a union of compact toric varieties sitting inside the hypertoric
$\MAt$, glued together along toric subvarieties as prescribed by the combinatorics
of the polytopes $P_U$ and their intersections in $\tddr$.  

\begin{example}\label{pictures}
Consider the two hyperplane arrangement pictured below, with all hyperplanes having primitive
normal vectors.  Note that there are two primitive vectors to choose from for each hyperplane
(one must choose a direction), but the corresponding hypertoric varieties and their cores will be independent of these choices by Proposition \ref{flip}.
\begin{figure}[h]
  \centerline{\epsfig{figure=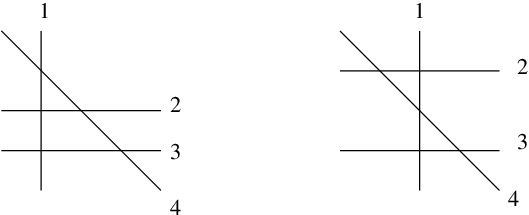, height=2cm}}
\end{figure}
\noindent
In the first picture, the core consists of a $\C P^2$ (the toric variety associated to a triangle)
and a $\C P^2$ blown up at a point (the toric variety associated to a trapezoid)
glued together along a $\C P^1$ (the toric variety associated to an interval). In the second
picture, it consists of two copies of $\C P^2$ glued together at a point.
\end{example}

\begin{remark}\label{cover}
Each of the core components $\mathcal E_U$ is a lagrangian subvariety of $\MAt$,
therefore its normal bundle in $\MAt$ is isomorphic to its cotangent bundle.  Furthermore,
each $\mathcal E_U$ has a $\Td$-invariant algebraic tubular neighborhood in $\MAt$
(necessarily isomorphic to the total space of $T^*\mathfrak X_U$), and these neighborhoods
cover $\MAt$.  Thus $\MAt$ is a union of cotangent bundles of toric varieties,
glued together equivariantly and symplectically in a manner prescribed by the combinatorics
of the bounded chambers of $\At$.
It is possible to take Proposition \ref{corecomps} and Equation
\eqref{pu} as a definition of $\mathfrak X_U$, and this remark as a definition of $\MAt$.
The affine variety $\MA$ may then be defined as the spectrum of the ring of global
functions on $\MAt$.
\end{remark}

\begin{remark}\label{why}
Though Propositions \ref{coreretract}, \ref{corecomps}, and \ref{wholecore} 
appear in the literature only for $\At$ simple, this hypothesis should not be necessary.
\end{remark}
\end{subsection}

\begin{subsection}{The Lawrence toric variety}\label{lawrence}
Let $$\mathfrak{B}(\At) = \cot\mod_{\!\!\a}\Tk.$$
This variety is a GIT quotient of a vector space by the linear action of a torus,
and is therefore a toric variety.  Toric varieties that arise in this way are
called {\bf Lawrence toric varieties}.  The following proposition is proven
in \cite[\S 6]{HSt}.

\begin{proposition}\label{lawrenceretract}
The inclusion 
$$\MAt = \mu^{-1}(0)\mod_{\!\!\a}\Tk\,\,\hookto\,\,\cot\mod_{\!\!\a}\Tk 
= \mathfrak{B}(\At)$$
is a $\Td_\R$-equivariant homotopy equivalence.
\end{proposition}

This Proposition is proven by showing that any toric variety retracts equivariantly
onto the union of those $\Td$-orbits whose closures are compact.  In the
case of the Lawrence toric variety, this is nothing but the core $\mathfrak L(\At)$.
\end{subsection}

\begin{subsection}{All the GIT quotients at once}\label{all}
Given $\a\in\tkdz$, we may define {\bf stable} and {\bf semistable} sets
$$\vast\subs\vass\subs\Cn$$ as in Section \ref{hkas},
and the toric variety $\X_\a = \Cn\mod_{\!\!\a}\Tk$
may be defined as the categorical quotient of $\vast$ by $\Tk$.
In analogy with Section \ref{hkas}, we will call $\a$ {\bf generic} if the $\a$-stable
and $\a$-semistable sets of $\Cn$ coincide.
In this case the categorical
quotient will be simply a geometric quotient, and $\X_\a$ will be the toric orbifold
corresponding to the polytope $P_\emptyset$ of Section \ref{core}.
We consider two characters to be equivalent if their stable sets are the same,
and note that there are only finitely many equivalence classes of characters, given by the various
combinatorial types of $P_\emptyset$ for different simplifications $\At$ of $\A$.
Let $\a_1,\ldots,\a_m$ be a complete list of representatives of equivalence classes
for which\footnote{Though $\mu^{-1}(\la)^{\a-st}$ is never empty,
$\vast$ sometimes is.} $\emptyset\neq\vast=\vass$.

Let $\vlf$ be the set of vectors in $\Cn$ on which $\Tk$ acts locally freely,
meaning with finite stabilizers.  For any character $\a$ of $\Tk$,
the stable set $(\Cn)^{\a-st}$ is, by definition, contained in $\vlf$.  Conversely,
every element of $\vlf$ is stable for some generic $\a$ \cite[1.1]{P4}, therefore
\begin{equation*}\label{stable}
\vlf = \bigcup_{i=1}^m\,\,(\Cn)^{\a_i-st}.
\end{equation*}
We define the {\em nonhausdorff} space
$$\xlf = \vlf/\,\Tk = \bigcup_{i=1}^m\,\,(\Cn)^{\a_i-st}/\Tk = \bigcup_{i=1}^m\,\,\X_{\a_i}$$
to be the union of the toric varieties $\X_{\a_i}$ along the open
loci of commonly stable points.

For an arbitrary $\la\in\tkd$, consider the projection
$$\pi_\la:\mu^{-1}(\la)\,\,\hookto\,\,\cot\to\Cn.$$
The following proposition is proven in \cite[1.3]{P2}.

\begin{proposition}\label{lf}
If $\la$ is a regular value of $\mu$, then $\pi_\la$ has image $\vlf$,
and the fibers of $\pi_\la$ are affine spaces of dimension $d$.
\end{proposition}

\begin{corollary}\label{wow}
The variety $\M_\la = \mu^{-1}(\la)/\Tk$ is an affine bundle
over $\xlf = \vlf/\Tk$.
\end{corollary}

It follows from Corollary \ref{wow} that the natural projection $\Ml\to\xlf$ is
a weak homotopy equivalence, meaning that it induces isomorphisms
on all homotopy and homology groups.  It is not a homotopy equivalence in the
ordinary sense because it does not have a homotopy inverse--in particular, it
does not admit a section.

\begin{example}
Consider the action of $\cs$ on $\C^2$ by the formula $t\cdot(z_1,z_2) = (tz_1,t^{-1}z_2)$.
A multiplicative character of $\cs$ is given by an integer $\a$, and that character will be generic
if and only if that integer is nonzero.  The equivalence class of generic characters
will be given by the sign of that integer, so we let $\a_1=-1$ and $\a_2=1$.
The corresponding stable sets
will be
$$(\C^2)^{\a_1-st} = \C^2\smallsetminus\{z_1=0\}\text{   and   }
(\C^2)^{\a_2-st} = \C^2\smallsetminus\{z_2=0\}.$$
The corresponding toric varieties $\X_{\a_1}$ and $\X_{\a_2}$ will both be isomorphic
to $\C$, and $\xlf = \X_{\a_1}\cup\X_{\a_2}$ will be the (nonhausdorff)
union of two copies of $\C$
glued together away from the origin.

The moment map $$\mu:\C^2\times(\C^2)^\vee\to\tkd\cong\C$$
is given in coordinates by $\mu(z,w) = z_1w_1 - z_2w_2$.
The hypertoric variety $\M_\a = \mu^{-1}(0)\mod\Tk$ at a generic character is isomorphic to 
$T^*\C P^1$, and its core is the zero section $\C P^1$.  
It is diffeomorphic to $\M_{\la} = \mu^{-1}(\la)/\cs$, which is,
by Corollary \ref{wow}, an affine bundle over $\xlf$.  If we trivialize this affine bundle
over the two copies of $\C$, we may write down a family of affine linear maps $\rho_z:\C\to\C$
such that, over a point $0\neq z\in\C$, the fibers of the two trivial bundles are glued
together using $\rho_z$.  Doing this calculation, we find that $\rho_z(w) = w + z^{-2}$.
\end{example}

\begin{remark}
Both Proposition \ref{coreretract} and Corollary \ref{wow} show that a hypertoric variety
is equivariantly (weakly) homotopy equivalent to a union of toric orbifolds.  
In the case of Proposition
\ref{coreretract} those toric orbifolds are always compact, and glued together along closed
toric subvarieties.  In the case of Corollary \ref{wow} those toric orbifolds may or may not
be compact, and are glued together along Zariski open subsets to create something
that has at worst orbifold singularities, but is not Hausdorff.  In general, there is no relationship
between the collection of toric varieties that appear in Proposition \ref{coreretract}
and those that appear in Corollary \ref{wow}.
\end{remark}

\begin{remark}
Corollary \ref{wow} generalizes to abelian quotients of cotangent bundles of arbitrary varieties,
rather than just vector spaces \cite[1.4]{P4}.
A more complicated statement
for nonabelian groups was used by Crawley-Boevey and Van den Bergh \cite{CBVdB}
to prove a conjecture of Kac about counting quiver representations over finite fields.
\end{remark}
\end{subsection}
\end{section}

\begin{section}{Cohomolgy}\label{cohomology}
In this Section we discuss the cohomology of 
the orbifold $\MAt$ and the intersection cohomology of the
singular variety $\MA$, focusing on the connection to the combinatorics of matroids.
In Section \ref{abel} we explain how hypertoric varieties can be used
to compute cohomology rings of nonabelian algebraic symplectic quotients,
as defined in Section \ref{hkas}.
There are a number of results on the cohomology of hypertoric varieties that
we won't discuss, including computations of 
the intersection form on the $L^2$-cohomology of $\MAt$ \cite{HSw} and
the Chen-Ruan orbifold cohomology ring of $\MAt$ \cite{GH,JT}.

\begin{subsection}{Combinatorial background}
A {\bf simplicial complex} $\D$ on the set $\otn$ is a collection of subsets of $\otn$,
called faces, such that a subset of a face is always a face.  Let $f_i(\D)$ denote the number
of faces of $\D$ of order $i$, and define the {\bf $\mathbf h$-polynomial} 
$$h_\D(q) := \sum_{i=0}^d f_i q^i (1-q)^{d-i},$$
where $d$ is the order of the largest face of $\D$.  Although the numbers $f_i(\D)$ 
are themselves very natural to consider, it is unclear from the definition above why we
want to encode them in this convoluted way.  The following equivalent construction of the
$h$-polynomial is less elementary but better motivated.

To any simplicial complex
one associates a natural graded algebra, called the {\bf Stanley-Reisner ring}, defined
as follows:
$$\SR(\D) := \C[e_1,\ldots,e_n]\Big/\<\prod_{i\in S}e_i\bigmid S\notin\D\>.$$
In order to agree with the cohomological interpretation that we will give to this
ring in Theorem \ref{smoothring}, we let the generators $e_i$ have degree 2.
Consider the Hilbert series 
$$\operatorname{Hilb}(\SR(\D), q) := \sum_{i=0}^\infty \dim \SR^{2i}(\D) q^i,$$ 
which may be expressed as a rational function in $q$.  
The following proposition (see \cite[\S II.2]{St})
says that the $h$-polynomial is the numerator of that rational function.

\begin{proposition}\label{num}
$\operatorname{Hilb}(\SR(\D), q) = h_\D(q)/(1-q)^d$.
\end{proposition}
\end{subsection}

\begin{subsection}{Cohomology of $\mathbf\MAt$}
Let $\D_\A$ be the simplicial complex consisting of all sets $S\subs\otn$ such that
the normal vectors $\{a_i\mid i\in S\}$ are linearly independent.
This simplicial complex is known as the {\bf matroid complex}
associated to $\A$.
The Betti numbers of $\MAt$ were computed in \cite[6.7]{BD}, but the following
combinatorial interpretation was first observed by \cite[1.2]{HSt}.
Let $$\Poin_{\MAt}(q) = \sum_{i=0}^d\dim H^{2i}(\MAt)\, q^i$$ be the even degree 
Poincar\'e polynomial of $\MAt$.

\begin{theorem}\label{smoothbetti}
The cohomology of $\MAt$ vanishes in odd degrees, and $$\Poin_{\MAt}(q) = h_{\DA}(q).$$
\end{theorem}

Theorem \ref{smoothbetti} is a consequence of the following stronger result.

\begin{theorem}\label{smoothring}
There is a natural isomorphism of graded rings $H^*_{\Td}(\MAt)\cong\SR(\DA)$.
\end{theorem}

The action of $\Td$ on $\MAt$ is equivariantly formal \cite[2.5]{K1}, therefore
the Hilbert series of $H^*_{\Td}(\MAt)$
is equal to $\Poin_{\MAt}(q)/(1-q)^d$, and Theorem \ref{smoothbetti} follows immediately
from Proposition \ref{num}.  Theorem \ref{smoothring} was proven for $\At$ smooth
in \cite[2.4]{K1} from the perspective of Section \ref{core}, and in the general case in 
\cite[1.1]{HSt}
from the perspective of Section \ref{lawrence}.  Here we give a new, very short proof, from
the perspective of Section \ref{all}.\\

\begin{proofsmoothring}
By Proposition \ref{diffeo}, Remark \ref{eq}, and Corollary \ref{wow}, 
$$H^*_{\Td}(\MAt) \cong H^*_{\Td}(\Ml) \cong H^*_{\Td}(\xlf) 
\cong H^*_{\Td}(\vlf/\Tk) \cong H^*_{\Tn}(\vlf).$$
Given a simplicial complex $\D$ on $\otn$, Buchstaber and Panov build a $\Tn$-space 
$\mathcal Z_\D$ called the moment angle complex with the property that
$H^*_{\Tn}(\mathcal Z_\D) \cong \SR(\D)$ \cite[7.12]{BP}.  In the case of the matroid complex $\DA$,
there is a $\Tn$-equivariant homotopy equivalence $\mathcal Z_{\DA}\simeq\vlf$ \cite[8.9]{BP},
which completes the proof.
\end{proofsmoothring}
\end{subsection}

\vspace{-\baselineskip}
\begin{subsection}{Intersection cohomology of $\mathbf \MA$}
The singular hypertoric variety  $\MA = \M_{0,0}$ is contractible, 
hence its ordinary cohomology is trivial.
Instead, we consider intersection cohomology, a variant of cohomology introduced
by Goresky and MacPherson which is better at probing the topology of singular varieties 
\cite{GM1,GM2}.  Let $$\Poin_{\MA}(q) = \sum_{i=0}^{d-1}\dim\IH^{2i}(\MA)\, q^i$$
be the even degree intersection cohomology Poincar\'e polynomial of $\MA$.  We will interpret
this polynomial combinatorially with a theorem analogous to Theorem \ref{smoothbetti}.

A minimal nonface of $\DA$ is called a {\bf circuit}.  
Given an ordering $\sigma$ of $\otn$, define a $\mathbf\sigma$-{\bf broken circuit}
to be a circuit minus its smallest element with respect to the ordering $\sigma$.
The {\bf$\mathbf\sigma$-broken circuit complex} $\bcsda$ is defined to be the collection
of subsets of $\otn$ that do not contain a $\sigma$-broken circuit.  Though the
simplicial complex $\bcsda$ depends on the choice of $\sigma$, its $h$-polynomial does not.
The following theorem was proved by arithmetic methods in \cite[\S 4]{PW}.

\begin{theorem}\label{singbetti}
The intersection cohomology 
of $\MA$ vanishes in odd degrees, and $$\Poin_{\MA}(q) = h_{\bcsda}(q).$$
\end{theorem}

Given the formal similarity of Theorems \ref{smoothbetti} and \ref{singbetti},
it is natural to ask if there is an analogue of Theorem \ref{smoothring} in the central case.
The most naive guess is that the equivariant cohomology $\IH^*_{\Td}(\MA)$
is naturally isomorphic to the Stanley-Reisner ring $\SR(\bcsda)$, but this guess is problematic
for two reasons.  The first is that intersection cohomology generally does not admit 
a ring structure, and therefore such an isomorphism would be surprising.  The second
and more important problem is that the ring $\SR(\bcsda)$ depends on $\sigma$,
while the vector space $\IH^*_{\Td}(\MA)$ does not.  
Since the various rings $\SR(\bcsda)$ for different
choices of $\sigma$ are not naturally isomorphic to each other, they cannot all be naturally
isomorphic to $\IH^*_{\Td}(\MA)$, even as vector spaces.  These problems can be addressed
and resolved by the following construction.

Let $R(\A) = \C[a_1^{-1},\ldots,a_n^{-1}]$ be the subring of the ring of all rational
functions on $\Cn$ generated by the inverses of the linear forms that define the hyperplanes
of $\A$.  There is a surjective map $\varphi$ from $\C[e_1,\ldots,e_n]$ to $R(\A)$ taking $e_i$
to $a_i^{-1}$.  Given a set $S\subs\otn$ and a linear relation of the form
$\sum_{i\in S}c_i a_i = 0$, the element 
$k_S=\sum_{i\in S} c_i \prod_{j\in S\smallsetminus \{i\}}e_j$ lies in the kernel $I(\A)$ of $\phi$,
and in fact $I(\A)$ is generated by such elements.  Since $k_S$ is clearly homogeneous,
$R(\A)$ is a graded ring, with the usual convention of $\deg e_i = 2$ for all $i$.
The following proposition, proven in \cite[4]{PS}, states that the ring $R(\A)$ is a simultaneous
deformation of the various Stanley-Reisner rings $\SR(\bcsda)$.

\begin{proposition}\label{ringdef}
The set $\{k_S\mid S\text{ a circuit }\}$ is a universal
Gr\"obner basis for $I(\A)$, and the choice of an ordering $\sigma$ of $\otn$
defines a flat degeneration of $R(\A)$ to the Stanley-Reisner ring $\SR(\bcsda)$.
\end{proposition}

\begin{example}
Let $d=2$, identify $\td_\R$ with $\R^2$, and let
$$a_1 = \binom{1}{0},\,\,\,\,\, a_2 = a_3 = \binom{0}{1},\,\,\, \text{ and } a_4 = \binom{-1}{-1}.$$
The two arrangements pictured in Example \ref{pictures} are two different simplifications
of the resulting central arrangement $\A$.
We then have
$$R(\A) \cong \C[e_1,\ldots,e_4]\,\big/
\<\, e_2-e_3,\,\,\, e_1e_2+e_1e_4+e_2e_4,\,\,\, e_1e_3+e_1e_4+e_3e_4\,\>.$$
By taking the initial ideal with respect to some term order, we get the Stanley-Reisner
ring of the corresponding broken circuit complex.
\end{example}

In Theorem \ref{singring}, proven in \cite{BrP}, 
we show that $R(\A)$ replaces the Stanley-Reisner ring in the ``correct"
analogue of Theorem \ref{smoothring}.

\begin{theorem}\label{singring}
Suppose that $\A$ is unimodular.  The equivariant intersection cohomology sheaf $\IC_{\Td}(\MA)$
admits canonically the structure of a ring object in the bounded equivariant
derived category of $\MA$.  This induces a ring structure on $\IH^*_{\Td}(\MA)$,
which is naturally isomorphic to $R(\A)$.
\end{theorem}

The problems of classifying $h$-polynomials of matroid complexes and their broken
circuit complexes remain completely open.  Hausel and Sturmfels explore the restrictions
on $h_{\DA}(q)$ imposed by Theorem \ref{smoothring} in \cite[\S 7]{HSt}, and Webster and
the author consider the combinatorial implications of applying the decomposition theorem
for perverse sheaves to the map $\MAt\to\MA$ \cite[\S 5]{PW}.
In both cases one obtains results which admit independent, purely
combinatorial proofs, but which are illuminated by their geometric interpretations.
\end{subsection}

\begin{subsection}{Abelianization}\label{abel}
As in Section \ref{hkas}, let $G$ be a reductive complex algebraic group acting linearly
on a complex vector space $V$, and let $T\subs G$ be a maximal torus.
We need the further technical assumption that $V$ has no nonconstant $T$-invariant functions,
which is equivalent to asking that any GIT quotient of $V$ by $T$
is projective.
The inclusion of $T$ into $G$ induces a surjection $\gd\surj\t^*$, which restricts to an
inclusion of $Z(\gd)$ into $\t^*$.  Thus a pair of parameters $(\a,\la)\in Z(\gd_\Z)\times Z(\gd)$
may be interpreted as parameters for $T$ as well as for $G$.
Suppose given $\a\in Z(\gd_\Z)$ such that $(\a,0)$ is generic for both $G$ and $T$,
so that the symplectic quotients
$$\Maz(G)\,\,\,\,\text{and}\,\,\,\,\Maz(T)$$ are both orbifolds.
Our first goal for this section is to describe the cohomology of $\Maz(G)$ in terms
of that of $\Maz(T)$.

Both $\Maz(G)$ and $\Maz(T)$ inherit actions of the group $\cs$ induced by scalar multiplication
on the fibers of the cotangent bundle of $V$.
Let $$\Phi(G):H^*_{G\times\cs}(T^*V)\to H^*_{\cs}(\Maz(G))$$
and
$$\Phi(T):H^*_{T\times\cs}(T^*V)\to H^*_{\cs}(\Maz(T))$$
be the {\bf equivariant Kirwan maps}, induced by the $\cs$-equivariant inclusions of
$\mu_G^{-1}(0)^{\a-st}$ and $\mu_T^{-1}(0)^{\a-st}$ into $T^*V$.
The map $\Phi(T)$ is known to be surjective \cite[4.5]{HP1},
and $\Phi(G)$ is conjectured to be so, as well.
The abelian Kirwan map $\Phi(T)$ makes the equivariant cohomology ring
$H^*_{\cs}(\Maz(T))$ into a module over
$H^*_{T\times\cs}(T^*V)$.
The Weyl group $W = N(T)/T$ acts both on the source and the target of $\Phi(T)$,
and the map is $W$-equivariant. 

Let $\D\subs\t^*$ be the set of roots of $G$ (not to be confused with the simplicial
complexes $\D$ that we discussed earlier), and consider the $W$-invariant class
$$e\,\, = \,\,\prod_{\beta\in\D}\beta \,(x-\beta) \,\,\in\,\, \Sym\t^*\otimes\C[x] 
\,\,\cong\,\, H^*_{T\times\cs}(T^*V).$$
The following theorem was proven in \cite[2.4]{HP}.

\begin{theorem}\label{ab}
If $\Phi(G)$ is surjective, then there is a natural isomorphism
$$H^*_{\cs}(\Maz(G)) 
\,\cong\, H^*_{\cs}(\Maz(T))^W\big/\operatorname{Ann}(e),$$
where $\operatorname{Ann}(e)$ is the ideal of classes annihilated by $e$.
\end{theorem}

We note that the abelian quotient $\Maz(T)$ is a hypertoric variety, and the $\cs$-equivariant
ring $H^*_{\cs}(\Maz(T))$ was explicitly described in \cite[4.5]{HP1} and \cite[3.5]{HH}.
Thus, modulo surjectivity of the Kirwan map, Theorem \ref{ab} tells us how to compute the cohomology ring of arbitrary symplectic quotients constructed in the manner of Section 
\ref{hkas}.  In \cite[\S 4]{HP}, this method was applied to compute the $\cs$-equivariant
cohomology rings of hyperpolygon spaces, a result which originally appeared in \cite[3.2]{HP2}
as an extension of the nonequivariant computation in \cite[7.1]{K2}.

Although the proof of Theorem \ref{ab} uses the $\cs$-action in a crucial way, Hausel
has conjectured a simpler, nonequivariant version.  Let $\Phi_0(G)$ be the map obtained
from $\Phi(G)$ by setting the equivariant parameter $x$ to zero, and let
$$e_0\,\, = \,\,\prod_{\beta\in\D}\beta\,\,\in\,\, \Sym\t^*
\,\,\cong\,\, H^*_{T}(T^*V).$$
Note that $e_0$ is {\em not} the class obtained from $e$ by setting $x$ to zero, rather
it is a square root of that class.

\begin{conjecture}\label{noneq}
If $\Phi_0$ is surjective, then there is a natural isomorphism
$$H^*(\Maz(G)) 
\,\cong\, H^*(\Maz(T))^W\big/\operatorname{Ann}(e_0).$$
\end{conjecture}

We end by combining Conjecture \ref{noneq} with Theorem \ref{singring}
to produce a conjecture that would put a ring structure on the intersection
cohomology groups of $\Mzz(G)$.
The hypothesis that $\A$ be unimodular in Theorem \ref{singring} is equivalent
to requiring that the orbifold resolution $\MAt$ of $\MA$ is actually smooth.
The analogous assumption in this context is that 
$\Maz(G)$ and $\Maz(T)$ are smooth for generic choice of $\a$.

\begin{conjecture}
Suppose that $\Maz(G)$ and $\Maz(T)$ are smooth for generic $(\a,0)$.  Then
The intersection cohomology sheaf $\IC(\Mzz(G))$
admits canonically the structure of a ring object in the bounded
derived category of $\Mzz(G)$, and
there is a natural ring isomorphism
$$\IH^*(\M_{0,0}(G)) 
\,\cong\, \IH^*(\Mzz(T))^W\big/\operatorname{Ann}(e_0).$$
\end{conjecture}
\end{subsection}
\end{section}

\footnotesize{

}
\end{document}